\documentclass{amsart}

\usepackage[latin1]{inputenc}
\usepackage{latexsym}
\usepackage{amsfonts}
\usepackage{amssymb}
\usepackage{color}
\usepackage{dsfont}
\usepackage{varioref}
\usepackage{epsfig}
\usepackage{hyperref}

\newtheorem{Lem}{Lemma}

\newtheorem{Thm}[Lem]{Theorem}

\newtheorem{Prop}[Lem]{Proposition}
\newtheorem{Cor}[Lem]{Corollary}

\newtheorem{Rem}[Lem]{Remark}

\newtheorem{Que}{Question}

\newcommand{\sptext}[1]{\ \text{#1} \ }

\DeclareMathOperator{\codim}{codim}

\newcommand{\es}{\textup{es}}

\newcommand{\NN}{{\mathds N}}

\newcommand{\CC}{{\mathds C}}
\newcommand{\RR}{{\mathds R}}

\newcommand{\smcdot}{{\textup{$\cdot$}}}
\newcommand{\smbullet}{{\textup{\tiny $^\bullet$}}}

\begin{document}

\title[Surfaces with Many Solitary Points]{Surfaces with Many Solitary
  Points}
\date{\today}
\author{Erwan Brugall\'e}
\address{Univ.\ Paris 6, 175 rue du Chevaleret, 75 013 Paris, France}
\email{brugalle@math.jussieu.fr}
\author{Oliver Labs}
\address{Universit\"at des Saarlandes, Mathematik und Informatik, Geb\"aude 2.4, D-66123 Saar\-br\"ucken, Germany}
\email{Labs@math.uni-sb.de, mail@OliverLabs.net}


\begin{abstract}
  It is classically known that a real cubic surface in $\RR P^3$ cannot have
  more than one solitary point (or $A_1^\smbullet$-singularity, locally
  given by $x^2+y^2+z^2=0$) whereas it can have up
  to four nodes (or $A_1^-$-singularities, locally given by $x^2+y^2-z^2=0$).
  We show that on any surface of degree $d\ge 3$ in $\RR P^3$ the maximum
  possible number of solitary points is strictly smaller than the maximum
  possible number of nodes.
  
  Conversely, we adapt a construction of Chmutov to obtain surfaces with many
  solitary points by using  a refined version of Brusotti's Theorem.
  Combining lower and upper bounds, we deduce: $\frac{1}{4}d^3 + o(d^3)\le
  \mu^3(A_1^\smbullet, d) \le \frac{5}{12}d^3  + o(d^3)$, where
  $\mu^3(A_1^\smbullet, d)$ denotes the maximum possible number of solitary
  points on a real surface of degree $d$ in $\RR P^3$.
  Finally, we adapt this construction to get real algebraic surfaces in $\RR
  P^3$ with many singular points of type $A_{2k-1}^\smbullet$ for all $k\ge 1$.
\end{abstract}

\maketitle

\section*{Introduction}

An \emph{ordinary double point}, or $A_1$-singularity,  of a hypersurface $f=0$
in $\RR P^n$ or $\CC 
P^n$ is a non-degenerate singular point $p$ of $f$; i.e.\ $f$ and all its
partial derivatives vanish at $p$, but
the
hessian matrix
$H_f (p)= ({\partial^2 f} \ / \ {\partial x_i\partial x_j}(p))_{i,j=0\dots n}$
is of rank $n$.
In $\RR P^3$, there are exactly two real types of ordinary double points:
we call the ones which can be given locally by the affine equation
$x^2+y^2-z^2=0$ \emph{nodes} or \emph{$A_1^-$-singularities},
and the others, locally given by $x^2+y^2+z^2=0$,
\emph{solitary ordinary double points},
\emph{$A_1^\smbullet$-singularities}, or \emph{solitary points} for short.

\cite{realArrNodes} showed by construction that for large degree $d$ the
currently known maximum number of complex singularities on a surface of degree
$d$ in $\CC P^3$ \cite{chmuP3} can also be achieved with a real surface with
only real singularities. 
All real singularities appearing in their construction are nodes.
In the present paper, we consider solitary points instead.

We denote the maximum possible number of complex $A_1$-singularities on a
complex hypersurface of degree $d$ in $\CC P^3$ by $\mu^3(A_1,d)$, and
similarly for the real $A_1^-$- and $A_1^\smbullet$-singularities on real
surfaces in $\RR P^3$: $\mu^3(A_1^-,d), \ \mu^3(A_1^\smbullet,d)$.

\begin{Que}
  It is clear that the maximum possible number of complex ordinary double
  points is at least as large as the corresponding real numbers: 
  \[\mu^3(A_1^-,d), \ 
  \mu^3(A_1^\smbullet,d) \ \le \ \mu^3(A_1,d).\] 
  Is any of these inequalities strict?
\end{Que}

\begin{Que}
  Classical results on cubic surfaces \cite{Schl63} and quartic surfaces 
  \cite{rohnSolPtsQuart} in $\RR P^3$ show that we have: 
  \[\mu^3(A_1^\smbullet, 3)=1<4=\mu^3(A_1^-, 3) \ \ \text{and} \ \ 
  \mu^3(A_1^\smbullet, 4)=10<16=\mu^3(A_1^-,4).\] 
  These results suggest that it might be more difficult to have many solitary
  points on surfaces than to have many nodes. 
  Is this true for all $d\ge 3$?
\end{Que}

In this article, we answer those questions involving solitary points
affirmatively in Theorem \ref{thm:A1bA1m}:   
\[
  \text{If} \ d\ge 3 \ \text{then} \ \mu^3(A_1^\smbullet,d) \ < \
  \mu^3(A_1^-,d), \ \mu^3(A_1,d).
\]
However the case of real nodes remains open in general
although it is clear that $\mu^3(A_1^-,d) \le \mu^3(A_1,d)$ for all $d$.
In fact, $\mu^3(A_1^-,d) = \mu^3(A_1,d)$ is only known for $d=1,2,\dots,6$. 

The currently known lower bound for $\mu^3(A_1^\smbullet,d)$ is still far
from the best known upper bound.
In the third section of this article, we improve the previously known maximum
number of solitary points on a surface of degree $d$ in $\RR P^3$ by adapting
a construction of Chmutov and by using
Brusotti's
Theorem. 
Altogether, we show for $d\in\NN$ by combining lower bound
(Theorem \ref{low1}) and upper bound
(Corollary \ref{cor:UpBA1b}): 
\[
    \frac{1}{4}d^3 + o(d^3)\le \mu^3(A_1^\smbullet,d) \le \frac{5}{12}
    d^3 + o(d^3).
\]

Together with the known cases in low degree, we get table \ref{tabA1}
which provides an overview of the known bounds for the maximum possible number
of both variants of the real ordinary double points.
In that table, the upper bounds for the case of $A_1^-$-singularities are
simply the complex ones most of which are due to Miyaoka \cite{miyp3}, the
asymptotic lower bound was found in \cite{realArrNodes}.

\begin{table}[h]
\begin{center}
  \begin{tabular}{|r||c|c|c|c|c|c|c|c||c|}
    \hline
    \rule{0pt}{1.1em}degree $d$ & $1$ & $2$ & $3$ & $4$ & $5$ & $6$ & $7$ &
    $8$ & large $d$\\
    \hline
    \hline
    \rule{0pt}{1.1em}$\mu^3(A_1^\smbullet,d)\ge$ & $0$ & $1$ & $1$ & $10$ & $12$ & $29$ &
    $45$ & $81$ & $\approx \frac{1}{4}d^3$\\[0.1em]
    \hline
    \rule{0pt}{1.1em}$\mu^3(A_1^\smbullet,d)\le$ & $0$ & $1$ & $1$ & $10$ & $24$ & $48$ &
    $83$ & $134$ & $\approx \frac{5}{12}d^3$\\[0.1em]
    \hline
    \hline
    \rule{0pt}{1.1em}$\mu^3(A_1^-,d)\ge$ & $0$ & $1$ & $4$ & $16$ & $31$ & $65$ & $99$ &
    $168$ & $\approx \frac{5}{12}d^3$\\[0.1em]
    \hline
    \rule{0pt}{1.1em}$\mu^3(A_1^-,d)\le$ & $0$ & $1$ & $4$ & $16$ & $31$ & $65$ & $104$ &
    $174$ & $\approx \frac{4}{9}d^3$ \\[0.1em]
    \hline
  \end{tabular}
\end{center}
\caption{An overview of the known bounds for the maximum possible number of both
variants of the real ordinary double points on surfaces of degree $d$ in $\RR
P^3$: solitary points and nodes.}
\label{tabA1}
\end{table}

More generally, an $A_j$-singularity of a complex surface in $\CC
P^3$ is a singular point locally given by the equation $x^{j+1} + y^2
+ z^2=0$. 
If $k\ge 2$ then there are three (two if $k=1$) real types of
$A_{2k-1}$-singularities, and we call the one given locally by the real
equation $x^{2k} + y^2 + z^2=0$ an \textit{$A_{2k-1}^\smbullet$-point }.
In section
\ref{lowerB 2k-1}, we explain how to adapt our  method to
construct real surfaces in $\RR P^3$ with many
$A_{2k-1}^\smbullet$-points. More precisely, we prove that (see
Proposition \ref{many Al}) for $k,d\ge 1$:
$$\frac{1}{8k-4}d^3 + o(d^3)\le \mu^3(A_{2k-1}^\smbullet,d)
\le  \frac{4k}{12k^2 - 3}d^3 + o(d^3).$$
The upper bound is again Miyaoka's bound on the number of complex
$A_{2k-1}$-points
of a complex surface of degree $d$ in $\CC P^3$.

\vspace{2ex}
\textbf{Acknowledgements:} We are grateful to Fr\'ed\'eric Bihan, Michel
Coste, Ilia Itenberg, and Jean Jacque Risler  for valuable and stimulating
discussions.

\section{Plane Curves with Solitary Points}
\label{sec:PlaneCurves}

In our results on real surfaces with solitary points we will use some facts
about real plane curves with solitary points.
So, we give a brief overview about this classical subject.
As in the case of $A_1$-singularities of surfaces mentioned in the
introduction, there are exactly two real types of ordinary double points on a
real plane curve, also denoted by $A_1^\smbullet$ resp.\ $A_1^-$.

\subsection{Nodes}
\label{sec:PlCNodes}

The value $\mu^2(A_1^-, d)$ of the maximum possible number of nodes on a real
plane curve of degree $d$ has been known for a long time:
  $$\mu^2(A_1^-, d) =\frac{d(d-1)}{2}.$$
The upper bound is a consequence of the genus formula, and a
generic configuration of $d$ lines shows that this upper bound is
sharp.
The genus formula also shows that this bound can only be achieved with
arrangements of $d$ real lines no three of which meet in a point.

There is a classical theorem, the Brusotti Theorem, which shows that we can
smooth each of the ordinary double points of a plane curve independently.
Applied to the $d$ generic lines in the plane mentioned above, we may deduce
that for any integer $r$ between
$0$ and $\frac{d(d-1)}{2} $, there is a real plane curve of degree
$d$ in $\RR P^2$ with exactly $r$ nodes as its only singularities.

Let us denote by $\mathcal C(d)$ (resp. $\RR \mathcal C(d)$) the space
of complex (resp. real) algebraic
curves of degree $d$ in $\CC P^2$ (resp. $\RR P^2$).
These are projective spaces of dimension $\frac{d(d+3)}{2}$.
Brusotti's result is the following:
\begin{Thm}[Brusotti Theorem, usual formulation]\label{usual Brusotti}
  Let $C$ be a real algebraic curve of degree $d$ in $\RR P^2$ with ordinary
  double points as its only singularities. For any of these singularities,
  choose a local deformation.
  Then it is possible to vary the curve $C$ in the space
  $\RR \mathcal C(d)$ in such a way that all previously chosen
  deformations are realized.
\end{Thm}

This is the form of the theorem which is usually given
because it can be applied very easily.
It is a straightforward corollary of the following result which will be more
convenient for our purposes.

\begin{Thm}[Brusotti Theorem, for a proof see e.g.\ \cite{BenRis}]\label{Brusotti}
  Let $C$ be a complex algebraic curve of degree $d$ in $\CC P^2$ with
  ordinary double points 
  $p_1,\ldots,p_k$ as
  its only singularities. Then there exists
  a small neighborhood $V_i$ of
  $p_i$ in $\CC P^2$ for each $i$, and a small  neighborhood $V$ of
  $C$ in $\mathcal C(d)$ such that the analytic sets
  $$S_i=\{\widetilde C\in V \ |\ \widetilde
  C \textrm{ is non-singular except at some point in $V_i$ where it has
    an $A_1$}\}$$
  are all non-singular and intersect transversely.
  Moreover, the tangent space of $S_i$ at $C$ is
  $\{\widetilde C \in \mathcal C(d)\ | \ p_i\in \widetilde C\}$.
\end{Thm}

\subsection{Solitary Points}
\label{sec:PlCSolPts}

For solitary points, things are a bit more complicated, and the exact
maximum number
of solitary points is only known since the 80's.

\begin{Prop}\label{up1}
  Let $d\in\NN$. Then:
  $$\mu^2(A_1^\smbullet, d)  \le \frac{(d-1)(d-2)}{2} +1.$$
\end{Prop}
\begin{proof}
  According to Harnack's Theorem, a non-singular real
  algebraic curve of degree $d$ in $\RR P^2$ has at most
  $\frac{(d-1)(d-2)}{2} +1$ connected components.
  Now the result follows
  from the Brusotti Theorem.
\end{proof}

In most cases, this upper bound can be refined using the Petrovskii
inequality (see \cite{Pet} or \cite{V3}):

\begin{Prop}\label{up2}
  If $2, 4 \ne d\in\NN$ then:
  $$\mu^2(A_1^\smbullet, d)  \le \frac{(d-1)(d-2)}{2}.$$
\end{Prop}
\begin{proof}
  This bound is trivial for curves of odd degree, as one
  component of the curve is not contractible in $\RR P^2$.
  The Petrovskii inequality
  for plane curves implies that
  if a curve of degree $d=2k$ has $\frac{(d-1)(d-2)}{2} + 1$ ovals, then
  at least one of them contains another oval if $k\ge 3$.
\end{proof}

The union of two complex conjugated
lines is a real conic with one $A_1^\smbullet$-point. If 
$P_1(x,y)=0$ and $P_2(x,y)=0$  are real equations of two real conics intersecting in four real points,
then the 
real quartic with equation $P_1^2(x,y)+P_2^2(x,y)=0$ has four
$A_1^\smbullet$-points. 
Hence, one has $\mu^2(A_1^\smbullet, 2)=1$ and
$\mu^2(A_1^\smbullet, 4)=4$ (i.e.\  Proposition \ref{up1}
is sharp in degree $2$ and $4$). 
The proof that the upper bound given in Proposition \ref{up2} is sharp for any
other degree has first been given by Viro
in the
80's. 
As Viro's original proof was not available to us, we sketch Kenyon
and Okounkov's \cite{KenyonOkounkov} here:

\begin{Thm}[Viro, see \cite{V0}, \cite{KenyonOkounkov}, or see \cite{Sh5} for a proof of 
 a more general
  case]\label{thm:curvesExact} 
  If $d\ne 2, 4$ then:
  $$\mu^2(A_1^\smbullet, d)  = \frac{(d-1)(d-2)}{2}.$$
\end{Thm}
\begin{proof}
 Let $\varepsilon$ be a primitive $d^{\text{th}}$ root of unity, and define the
 polynomial 
 $\widetilde P_d(x,y)=\prod_{i,j=1}^d (\varepsilon^ix + \varepsilon^j
 y +1)$. Then, $\widetilde P_d(x,y)=P_d(x^d,y^d)$ where $P_d(x,y)$ is
 a real polynomial of degree $d$, and the curve with
equation $P_d(x,y)=0$ has exactly $\frac{(d-1)(d-2)}{2}$ real
 solitary non-degenerate double points.
\end{proof}

As in the case of nodes, it follows from Brusotti's Theorem that for any 
integer
$r$ between 0 and $\frac{(d-1)(d-2)}{2} $, there exists a real algebraic
plane curve of degree $d$ in $\RR P^2$ with exactly $r$ solitary nodes as 
 its
only singularities.

For what follows, we need to introduce a distinction among solitary
points of a real algebraic curve $P(x,y)=0$ of even degree in $\mathbb C^2$ :
those who are local minima of 
the function $(x,y)\mapsto P(x,y)$ and those who are local maxima.

\begin{Prop}[\cite{KenyonOkounkov},\cite{Mik11}]\label{harnack rat}
  Let $d\ge 2$ be even and let $P_d(x,y)=0$ be the polynomial constructed  in
  the  proof of Proposition
  \ref{thm:curvesExact}. Then
  $\frac{3d(d-2)}{8}$ solitary points of the curve $P(x,y)=0$ are local
  maxima of the polynomial $P(x,y)$ and 
  the  $\frac{(d-2)(d-4)}{8}$ other solitary points of the curve
  $P(x,y)=0$  are local minima of the polynomial $P(x,y)$.
\end{Prop}
\begin{proof}
  Kenyon's and Okounkov's polynomial $P_d(x,y)$ defines a Harnack curve, and
  Mikhalkin showed that these curves have the desired property.
\end{proof}

\section{Surfaces in $\RR P^3$: Upper Bounds}

In order to prove the upper bounds mentioned in the introduction, we need 
 ---
in analogy to the Brusotti Theorem in the case of plane curves --- a
result about smoothings of algebraic varieties.

By a \emph{smoothing} of a singular (real) algebraic hypersurface $X$ of
degree $d$ in $\CC P^n$, we mean a small perturbation of the
coefficients of $X$ such that the result is a non-singular (real)
algebraic hypersurface of degree $d$ in $\CC P^n$.
The Coste-Hironaka Theorem now says that one can always smooth a real projective
hypersurface  in such a way that no connected component disappears into the
complex world:

\begin{Thm}[Coste-Hironaka, \cite{Cos1}]\label{thm:Coste-Hironaka}
  Let $X$ be a singular real algebraic hypersurface in $\RR P^n$. Then
  there is a smoothing $\widetilde X$ of $X$ such that
  $$b_0(X)\le b_0(\widetilde X),$$
  where $b_0$ denotes the $0^\text{th}$ Betti number, i.e.\ the number of
  connected components.
\end{Thm}

\begin{Rem}
  In the special case of hypersurfaces with only solitary (ordinary double!)
  points as singularities, it is easy to prove this result using the
  construction given in the paper \cite{Cos1}:
  indeed, let
  $P(X_0,\ldots,X_n)=0$ be the equation of such a hypersurface in $\RR P^n$
  which does not have a singularity in the point $(1:0:\cdots:0)$ (which we
  may assume after a suitable change of coordinates).
  Then  define
  $$\widetilde P(X_0,\ldots,X_n) := P(X_0,\ldots,X_n) +  \sum_{i=1}^n
  \varepsilon_i X_i \frac{\partial P}{\partial
    X_i}(X_0,\ldots,X_n).$$ 
  Each singular point $p$ of $P=0$ will still be a point on
  $\widetilde{P}=0$. 
  Moreover, a short computation shows that there are $\varepsilon_i$ small
  enough such that $\widetilde{P}=0$ is non-singular in $p$ because of the
  hessian criterion for $A_1$-singularities.
  But this means that if the $\varepsilon_i$ are small enough then near each
  solitary point $p$, the hypersurface $P=0$ is
  smoothed into a small connected component of $\widetilde P=0$ homeomorphic
  to an $n$-sphere and containing $p$. 
\end{Rem}

As we know the homology of projective non-singular complex algebraic
hypersurfaces,
the Coste-Hironaka Theorem \ref{thm:Coste-Hironaka} combined with Smith Theory
(see \cite{Bred1}) implies the following corollary:

\begin{Cor}[Coste-Hironaka, \cite{Cos1}]
Let $X$ be a (possibly singular) real algebraic hypersurface of degree $d$ in
$\RR P^n$. Then 
$$b_0(X)\le  \frac{1}{2}\left(\frac{(d+1)^{n+1} -(-1)^{n+1}}{d} +n - (-1)^n
\right).$$ 
\end{Cor}

In the case of projective surfaces in $\RR P^3$ one can
improve the upper bound on the number of connected components
thanks to the Petrovskii-Oleinik inequality (see, e.g.,
\cite{degKhRokhlin}):

\begin{Cor}\label{cor:b0upper}
  Let $S$ be a (possibly singular) real algebraic surface of degree $d$ in $\RR P^3$. Then
  $$b_0(S)\le  \frac{5d^3 - 18d^2 +  25d}{12}.$$
\end{Cor}
\begin{Rem}
 Note that starting from
degree 5, the maximal possible value of $b_0(S)$ when $S$ is a
 real algebraic surface of degree $d$ is still unknown.
\end{Rem}

Applying Corollary \ref{cor:b0upper} in the special case of 
real surfaces with solitary double points in $\RR P^3$, we get:

\begin{Cor}\label{cor:UpBA1b}
  For $d\in\NN$, we have:
  $$\mu^3(A_1^\smbullet,d)\le
  \begin{array}{l@{\qquad}l}
    \left\lfloor\frac{5d^3 - 18d^2 +  25d}{12}\right\rfloor, & d \text{\ even,}\\[0.75em]
    \left\lfloor\frac{5d^3 - 18d^2 +
        25d}{12}\right\rfloor - 1,&  d \text{\ odd.}
  \end{array}$$
\end{Cor}
\begin{proof}
  In odd degree, we can subtract one because in that case at least one of the
  connected components from Corollary \ref{cor:b0upper} is not homeomorphic to a
  sphere.  
\end{proof}

Comparing this upper bound with the lower bound obtained in the case of nodes
(see Theorem 2 in \cite{realArrNodes} for a detailed formula) which is given by 
$$\mu^3(A_1^-,d)\ge
\begin{array}{l@{\qquad}l}
  \frac{5}{12}d^3 - \frac{13}{12}d^2 + o(d^2), & d \sptext{even,}\\[0.3em]
  \frac{5}{12}d^3 - \frac{14}{12}d^2 + o(d^2), & d \sptext{odd,}
\end{array}$$
we may deduce that one cannot reach the maximum number of nodes with surfaces
having only solitary points:
\begin{Thm}\label{thm:A1bA1m}
  The maximum possible number of solitary points on a surface of degree $d,
  d\ge3,$ in $\RR P^3$ is strictly smaller 
than the maximum possible number of nodes:
  \[\mu^3(A_1^\smbullet,d) \ < \ \mu^3(A_1^-,d), \ \mu^3(A_1,d).\]
\end{Thm}

We already mentioned in the introduction that this result is not very
surprising because it has been known for degree $3$ and $4$ for a long time. 
However, notice that the corresponding statement in the case of plane curves
does not hold: the maximum number of nodes on an irreducible curve of degree
$d$ equals the maximum number of solitary points on an irreducible curve of
degree $d$: in both cases, it is the genus of a smooth plane curve of degree
$d$, as mentioned earlier.

\section{Surfaces in $\RR P^3$: Lower Bounds by Constructions}
\label{sec:LowerBounds}

In the preceding section we have shown that the maximum possible number of
solitary points on a surface in $\RR P^3$ is less than the corresponding
number of nodes.
Here we improve the currently known maximum number of
solitary points.
Indeed, in this section we show:

\begin{Thm}\label{low1}
  Let $d\in\NN$. Then:
  $$\begin{array}{ll}
    \mu^3(A_1^\smbullet,d) \ge \frac{ (d-2)(2d^2 - 3d + 4)}{8} &\text{ if d is even,}\\[0.4em]
    \mu^3(A_1^\smbullet,d) \ge \frac{(d-1)^2(d-2)}{4} &\text{ if d is odd}.
  \end{array}$$
\end{Thm}

We prove Theorem \ref{low1} in section \ref{prflow1}.
It is based on Chmutov's method to construct singular complex
surfaces.
We thus explain this method first.
Then we discuss how to adapt it to obtain real algebraic surfaces with
solitary nodes; 
finally, we show in section \ref{optimality} that our result is asymptotically
the best that one can achieve using Chmutov's method.

\subsection{Known Constructions}
\label{sec:KnownConstr}

Notice that the previously best known lower bound for the maximum possible
number $\mu^3(A_1^\smbullet,d)$ of solitary points on a surface in $\RR P^3$
is far below $\frac14 d^3$.

Not many constructions are known.
Certainly the most sophisticated one is Shustin's variant of Viro's
patchworking method for the singular case.
This method yields
the  optimal result in the case of plane curves,
but already for complex surfaces in $\CC P^3$ it   
only yields
$\mu^3(A_1,d) \ge \frac1{6}d^3 + o(d^3)$ (see
\cite{shuWesPresc}).

There is another construction which is natural to consider:
we take a polynomial  $P_d(x,y)$ of degree $d$  in two
variables and we set
\[f(x,y,z) = P_d(x,y) + g(z),\]
where
$g(z)$ is a polynomial function of degree $d$
in  one variable $z$
with the maximum possible number $\lceil \frac{d}{2} \rceil$ of local
maxima
$z_i$ with value $g(z_i)=0$.
An \emph{even solitary point} of an affine plane curve
given by the equation $P(x,y)=0$ is a solitary point $(x,y)$ of
the curve $P(x,y)=0$ which is a local minimum of the polynomial $P(x,y)$,
i.e.\ locally at $p$ 
the graph of $P(x,y)$ looks like $z=x^2+y^2$.
We denote by $\es(P_d)$ the number of even solitary points
of the curve $P_d(x,y)=0$.
With these preliminaries, it is clear that the surface $f(x,y,z)=0$
has
$\lceil d/2\rceil \cdot \textup{es}(P_d)$
solitary points:
for each even solitary point $(a,b)$ of the affine plane curve $P_d(x,y)=0$,
we thus get a point $(a,b,z_i)$ of $f(x,y,z)=0$ which is locally of the form
$x^2+y^2+z^2$.
However, it is well known that for a curve of degree $d$, one has (see
\cite{V3}) 
$$\es(P_d)\le \frac{7}{16}d^2 + o(d^2), $$
so one cannot expect to construct in this way surfaces of degree $d$ with more than
$\frac{7}{32}d^3 + o(d^3) $ solitary points.
Combining this method with Chmutov's method we improve the
leading coefficient $\frac{7}{32}$.

\subsection{Chmutov's method}

We describe briefly how Chmutov constructed surfaces with many (complex)
ordinary double points in the 90's \cite{chmuP3}.
It is similar to the idea mentioned in the previous paragraph.
Despite its simplicity, Chmutov's surfaces still yield the best known lower bound
for the maximum number of ordinary double points on a complex surface of
degree $d\ge 13$.
The best known lower bound in the case of real nodes ($A_1^-$-singularities)
which we mentioned above and which equals the current lower bound in the
complex case is an adaption of Chmutov's construction to real nodes
\cite{realArrNodes}.
So, it is quite natural to try to adapt the method to solitary points.
However, we will see that this process is not completely straightforward,
and we will need a refined version of Brusotti's Theorem to make it work.

\subsubsection{Chmutov's Constructions}

Let $T_d(z)\in\RR[z]$ be the Tchebychev polynomial of degree $d$ with
$\lceil \frac{d-1}{2} \rceil$ extremal points with value $-1$
and $\lfloor \frac{d-1}{2} \rfloor$ with value $+1$.
This can either be defined recursively by $T_0(z) := 1$, \ $T_1(z) :=z$, \
$T_d(z) := 2\smcdot z\smcdot T_{d-1}(z) - T_{d-2}(z)$ for $d\ge2$, or
implicitly by $T_d(\cos(z)) = \cos(dz)$.
In \cite{chmuP3}, Chmutov used the Tchebychev polynomials to construct
surfaces in $\CC P^3$ with
$\approx \frac{5}{12}d^3$ (complex) nodes using the so-called \emph{folding
polynomials} $F^{A_2}_{d}(x,y)\in\RR[x,y]$ associated to the root-system
$A_2$:
\[\textup{Chm}_{d}^{A_2}(x,y,z) := F^{A_2}_{d}(x,y) + \frac{1}{2}(T_d(z) + 1).\]
The polynomials $F^{A_2}_{d}(x,y)$ have critical points with only three different critical
values: $0$, $-1$, and $8$.
The surface $\textup{Chm}_{d}^{A_2}(x,y,z)=0$ is singular exactly at those
points at which the critical values of $F^{A_2}_{d}(x,y)$ and
$\frac{1}{2}(T_d(z) + 1)$ sum up to zero (i.e., either both are $0$ or the
first one is $-1$ and the second one is $+1$).

\subsubsection{Adaption to Real Nodes}

In \cite{realArrNodes}, this construction was modified  to yield real 
surfaces $\textup{Chm}^{A_2}_{\RR,d}(x,y,z)=0$ with real nodes as singularities
by using the so-called \emph{real folding polynomials}
\[F^{A_2}_{\RR,d}(x,y) := F^{A_2}_{d}(x+iy, \ x-iy),\] where $i$ is the
imaginary number. 
It is not difficult  to see that the singularities of the surface
$\textup{Chm}^{A_2}_{\RR,d}(x,y,z)=0$ are indeed nodes, i.e.\ of type $A_1^-$,
by using the fact that the plane curve $F^{A_2}_{\RR,d}(x,y)=0$ is actually a
product of $d$ real lines no three of which meet in a common point.

\subsubsection{Adaption to Solitary Points}\label{sol points}

From the explanations in the previous paragraphs it is clear how to adapt
Chmutov's construction to yield solitary points:
we need to show the existence of a real polynomial $f(x,y)$ with many
local minima with
value $+1$ and local maxima with value $-1$.
More precisely:

\begin{Prop}\label{loc extr}
  Let $f(x,y)$ be a real polynomial of even (resp.\ odd) degree $d$ with
  $\alpha$ local minima with value $+1$,
  and $\beta$ local maxima with value $-1$.
  Then the affine surface defined by $f(x,y)-T_d(z)=0$ has
  $\frac12(\alpha\cdot d + \beta \cdot (d-2))$
  (resp.\ $\frac12 (\alpha+\beta)\cdot (d-1)$) solitary points. 
  The corresponding projective surface in $\RR P^3$ has at most $O(d^2)$
  additional  
singularities.
\end{Prop}

Notice that we cannot use a product of real lines such as
$F^{A_2}_{\RR,d}(x,y)$ as the polynomial $f(x,y)$ in order to obtain many
solitary points because it has the wrong critical values: the minima have
critical value $-1$ and the maxima $+1$.

\subsection{Proof of the Lower Bound of Theorem \ref{low1}}\label{prflow1}

We are now ready to prove Theorem \ref{low1} on the lower bound for
$\mu^3(A_1^\smbullet, d)$.
According to Chmutov's construction and in particular Proposition \ref{loc
  extr}, we have to construct polynomials in 2 variables whose graphs 
have many minima with 
value $+1$ and many maxima with value $-1$.
The existence of such polynomials is established by
Proposition
\ref{prop:RefinedBrusotti} below applied to polynomials constructed in
Propositions \ref{thm:curvesExact} and \ref{harnack rat}.
This completes the proof of Theorem \ref{low1}.\hfill\qed

\begin{Prop}\label{prop:RefinedBrusotti}
  Let $d\in\NN$ and $P(x,y)$ be a real polynomial of degree $d$ with $\alpha$
  (resp.  $\beta$)
  local minima  (resp. local maxima) with value 0. Then, there exists
  a real polynomial in two variables  of degree $d$
  with $\alpha$ local minima with value $+1$ and $\beta$
  local maxima with value $-1$.
\end{Prop}
\begin{proof}
  We start with the following observation:
  if $f(x,y)$ is a polynomial of degree $d$, then the graph of $f$,
  defined by the equation $f(x,y)-z=0$, is a special line in the space 
  $\mathcal C(d)$ of
  plane curves of degree $d$.
  Indeed, if $f(x,y)=\sum a_{i,j}x^iy^j$,
  then the section of the graph of $f$ by the hyperplane $z=t$ is
  given by the equation  $\sum a_{i,j}x^iy^j - t=0$.
  If $(a_{0,0}: a_{1,0}: a_{0,1}: \cdots :a_{0,d})$ are the coordinates in
  the space $\mathcal C(d)$, then the
  graph of $f$ can be
  parameterized by the line $t\mapsto (a_{0,0} - t:
  a_{1,0}: a_{0,1}: \cdots:a_{0,d})$. For $t=\infty$, this line passes
  through the point $(1:0:0:\cdots:0)$ which represents the multiple line
  $z^d$.
  Conversely, any  line in the space of plane curves of degree $d$
  passing through the point  $(1:0:0:\cdots:0)$ admits a
  parameterization of the form
  $t\mapsto (a_{0,0} - t:
  a_{1,0}: a_{0,1}: \cdots:a_{0,d})$ which 
  defines a polynomial
  $f(x,y)=\sum a_{i,j}x^iy^j$ of degree $d$.

Let us go back to the polynomial $P(x,y)$ of the Theorem. By
assumption, the curve defined by $P$ has $\alpha+\beta$ solitary points.  
Now we show that we can perturb the polynomial
  $P(x,y)$ in such a way
  that all local minima (resp.\ maxima) stay on the same level $a$
  (resp.\ $b$) with $a>b$ (see Figure \ref{pert}).
 
  \begin{figure}[h]
    \centering
    \begin{tabular}{c}
      \includegraphics[height=9cm, angle=0]{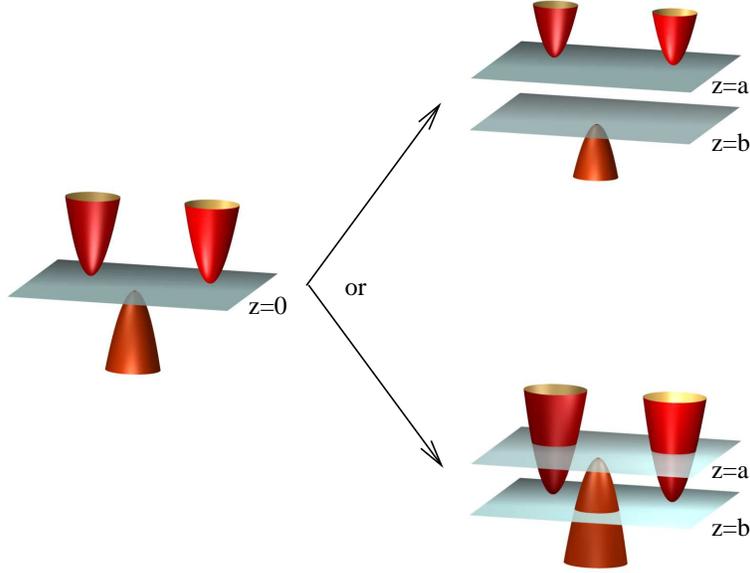}
    \end{tabular}
    \caption{Two ways to perturb $P(x,y)$.}
    \label{pert}
  \end{figure}

  For any  solitary point $p$ of the curve $P(x,y)=0$, we choose a
  small neighborhood $V(p)$ of 
  $p$ in $\RR P^2$ such that $V(p)\cap V(q)=\emptyset$ if $q\ne p$ is
  another 
  solitary point.
  We denote by $M(P)$ (resp.\ $m(P)$) the set of solitary points of
  the curve  $P(x,y)=0$
  corresponding to local maxima (resp.\ minima) of $P(x,y)$.
  Moreover, we denote by
  $\Sigma_{M(P)}$ (resp.\ $\Sigma_{m(P)}$) the stratum of real algebraic
  plane curves in $\mathcal C(d)$ in a small
  neighborhood of $P(x,y)=0$ with one solitary point in $V(p)$ for any
  $p\in 
  M(P)$ (resp.\ $m(P)$). Then, according to Brusotti's Theorem,
  $\Sigma_{M(P)}$  and  $\Sigma_{m(P)}$ are smooth and intersect
  transversely at the curve $P(x,y)=0$.
  Moreover, we have:
  $$\codim (\Sigma_{M(P)})= \beta, \ \  \ \codim (\Sigma_{m(P)})= \alpha$$
  and  
  $$ \codim (\Sigma_{M(P)}\cap\Sigma_{m(P)} )= \beta +\alpha \le
  \frac{(d-1)(d-2)}{2}=\frac{d(d+3)}{2}-(3d-1).$$
 
  One can suppose that $\alpha>0$ and $\beta>0$  otherwise the
  proposition is trivial. We denote by $L$ the line in the space
  $\mathcal C(d)$  passing through the curve
  $P(x,y)=0$ and $z^d$. By a simple
  dimension computation, we prove that we can perturb $L$ to a line
  $\widetilde L$
  still passing through $z^d=0$ and intersecting the stratum
  $\Sigma_{M(P)}$  and  $\Sigma_{m(P)}$ one after the other:
 
  Define the projection
  $$\begin{array}{cccc}
    \pi  : & \RR P^{\frac{d(d+3)}{2}}& \longrightarrow  & \RR P^{\frac{d(d+3)}{2}-1}
    \\  & (a_{0,0}:a_{1,0} : a_{0,1} : \cdots)& \longmapsto &(a_{1,0} : a_{0,1}
    : \cdots).
  \end{array}$$  
  None of the tangent
  spaces of $\Sigma_{M(P)}$  and  $\Sigma_{m(P)}$ contains the point
  $z^d=0$, so $ \pi\big(\Sigma_{M(P)}\big)$ and
  $\pi\big(\Sigma_{m(P)}\big)$ are non-singular and intersect transversely.
  Hence, we have:
  $$\codim (\pi\big(\Sigma_{M(P)}\big))= \beta - 1, \ \ \
  \codim (\pi\big(\Sigma_{m(P)}\big))= \alpha - 1,  $$
  $$  \codim (\pi\big(\Sigma_{M(P)}\cap\Sigma_{m(P)}\big) )= \beta +
  \alpha - 1$$
  and
  $$ \codim (\pi\big(\Sigma_{M(P)}\big) \cap\pi\big( \Sigma_{m(P)}\big)
  )= \beta +
  \alpha - 2.$$
  So, we have one degree of freedom to move from $\pi(P(x,y)=0)$ out of
  $\pi\big(\Sigma_{M(P)}\cap\Sigma_{m(P)}\big)$ staying in
  $\pi\big(\Sigma_{M(P)}\big) \cap\pi\big( \Sigma_{m(P)}\big)$ which
  means exactly that we can perturb $L$ to a line
  $\widetilde L$
  still passing through $z^d=0$ and intersecting the stratum
  $\Sigma_{M(P)}$  and  $\Sigma_{m(P)}$ one after the other.

  As
  $\pi\big(\Sigma_{M(P)}\big) \cap\pi\big( \Sigma_{m(P)}\big)\setminus
  \pi\big(\Sigma_{M(P)}\cap\Sigma_{m(P)}\big)$ has two connected
  components, we have two possible choices 
  to perturb $L$. One will correspond to move up (resp.\ down) the
  local maxima  (resp.\ minima) and the other will correspond to
  move up (resp.\ down) the
  local minima (resp.\ maxima), see Figure \ref{pert}. Choosing the
  latter possibility, we
  prove the proposition.
\end{proof}

The proposition can be interpreted as a refined version of Brusotti's
Theorem in a special case.
Indeed, it does not only show that we can perturb each solitary point of a
real plane curve $P(x,y)=0$ into one of the two topological possibilities, but it proves
that we can in addition put all solitary 
points which are deformed in the same topological way on the same
level of $P(x,y)$, 
i.e.\ transform the points into extremal points of the graph of $P(x,y)$ with the
same value.

\subsection{Optimality of our Construction}\label{optimality}

We now show that using Chmutov's method it is asymptotically not possible to
improve our lower bound obtained in Theorem \ref{low1}.
Let us denote by $\mu_\text{\it Ch}(d)$ the maximal possible number of
solitary points of a real algebraic surface of degree $d$ in $\RR P^3$
constructed using Chmutov's method.

\begin{Prop}\label{best}
  Let $d\in\NN$. Then: $$\mu_\text{\it Ch}(d) = \frac{1}{4}d^3+ o(d^3).$$
\end{Prop}
\begin{proof} 
  The result is an immediate corollary of Theorem \ref{low1},
  Proposition \ref{loc extr} and of the following Proposition
  \ref{prop:muExtr}.
\end{proof}

 Let us denote by $\mu_\text{\it extr}(d)$ the
  maximum possible number of local extrema of a real polynomial
  $f(x,y)$ of degree $d$.
We believe that the bound we establish now is known, but as we did not find a
reference for it, we include a proof here:

\begin{Prop}\label{prop:muExtr}
  With the notation of the preceding proof, we have:
  \[\mu_\text{\it extr}(d) \le \frac{1}{2}d^2+ o(d^2).\]
\end{Prop}
\begin{proof}
  Denote by $h$ the height function $(x,y,z)\mapsto z$
  defined on $\RR^3$.
  Let $f(x,y)$ be a real polynomial of degree $d$ and denote by
  $\nu_0(f)$ (resp.\  $\nu_1(f)$) the number 
  of local extrema (resp.\ hyperbolic critical points) of
  $f$. Consider a very large ball $B$ in $\mathbb R^2$
  containing all critical points of $f$, and consider $D(f)$ the
  intersection of the graph of $f$ with the cylinder with base
  $B$. Then, one can glue in $\RR^3$ a disk to $D(f)$ along its border
  $\partial D(f)$ adding  a
  number of critical points for $h$ which is at most linear in $d$. Then, we
  obtain a sphere $S^2$ and $h$ defines a Morse function on it. Hence,
  we have $\nu_0(f) - \nu_1(f) \le 2 + a d$ with $a$ some integer number.
 
  On the other hand, the number of real critical points of $f$ is not
  more than its number of complex critical points, which is equal to
  $(d-1)^2$.
  Taking all this together, we get:
  $\mu_\text{\it extr}(d) \le \frac{1}{2}d^2+ o(d^2).$
\end{proof}

\section{Higher singularities}\label{lowerB 2k-1}

Proposition \ref{loc extr} can also be applied to construct real algebraic
surfaces in $\RR P^3$ with many $A_{2k-1}^\smbullet$ singularities. 
The method is exactly the same as in 
section \ref{sol points}, but instead of Tchebychev polynomials,
we use polynomials with very degenerate critical points of critical
values $\pm1$.
The existence of such polynomials is guaranteed by applying the real
version of Dessins d'Enfants (e.g. see \cite{erwanMaxEvenOvals}) to the construction
in \cite{labsVarChmu}.
\begin{Lem}\label{degenerate Tcheb}
  Let $d,k\ge 1$.
  Then there is a real polynomial $T_d^{2k}(z)$ of degree $d$
  with $ \bigl[ \frac{d-1}{4k-2} \bigl] $ local maxima (resp.\ minima)
  which are 
  critical points of multiplicity $2k-1$ (resp. non-degenerate critical
  points) and with value $+1$
  (resp.\ $-1$).
\end{Lem}

\begin{Prop}\label{many Al}
  Let $k,d\ge 1$.  
  We have: 
  $$\frac{1}{8k-4} d^3 + o(d^3)  \le
  \mu^3(A_{2k-1}^\smbullet,d)\le  \frac{4k}{12k^2 - 3}d^3 + o(d^3)  .$$ 
\end{Prop}
\begin{proof}
  The upper bound is Miyaoka's bound. Let $f_d(x,y)$ be a real polynomial
  of 
  degree $d$ with $\alpha$ local minima with
  value $1$,
  with $\beta$ local maxima with value $-1$, and such that $\alpha+\beta=\frac{(d-1)(d-2)}{2}$. According to 
  Theorem \ref{thm:curvesExact} and
Proposition \ref{prop:RefinedBrusotti}, such a
polynomial exists.
The lower bound in the theorem is given by
  considering the surface with equation $f_d(x,y)-T_d^{2k}(z)=0$.
\end{proof}

\begin{Rem}
Using the method ``$P_d(x,y)+g(z)$=0'' described in section
\ref{sec:KnownConstr}, one could expect to obtain better lower bounds
for  $\mu^3(A_{2k-1}^\smbullet,d)$ as soon as $k\ge 2$. In this case,
$P_d(x,y)=0$ should be a plane curve with many even $A_{2k-1}^\smbullet$-points.
However, up to our knowledge, the currently known constructions give only
\begin{equation}\label{low2}
\mu^2(A_{2k-1}^\smbullet,d)\ge \frac{1}{4k}d^2
\end{equation}
which provide  lower bounds a bit worse than ours for
$\mu^3(A_{2k-1}^\smbullet,d)$. The lower bound (\ref{low2}) can be obtained by
considering the polynomials  $ T_d(x) -  \widetilde
T_d^{2k}(y)$ where $T_d(x)$ is the Tchebychev polynomial of degree $d$
and $\widetilde T_d^{2k}(y)$ is a polynomial of degree $d$ which has 
$ \lfloor \frac{d}{2k} \rfloor $
local minima which are critical
points of multiplicity $2k-1$ and with value $+1$. The existence of the
polynomials $\widetilde T_d^{2k}(y)$ can be proved with the same
technique as in Lemma \ref{degenerate Tcheb}.

In \cite[Proposition 3.5]{westenbergerManySimple}, Westenberger  claims that
$\mu^2(A_{2k-1}^\smbullet,d)\ge \frac{1}{4k-2}d^2.$
However, his proof of this proposition uses
\cite[Lemma 3.1]{westenbergerManySimple} which is wrong
for solitary points. 
Indeed, this lemma states that there exists a real algebraic curve with Newton
polygon the quadrangle with vertices $(0,1)$, $(0,2)$, $(1,0)$ and
$(2k-1,1)$, and with one $A_{2k-1}^\smbullet$ point.
However, such a curve cannot exist due to the following proposition.
The case of $A_1^\smbullet$-singularities is easy to verify by hand;
for the general statement, we need some more work: 
\begin{Prop}
For any $k\ge 1$, no real algebraic curve with Newton
polygon the quadrangle with vertices $(0,1)$, $(0,2)$, $(1,0)$ and
$(2k-1,1)$ can have an $A_{2k-1}^\smbullet$ point.
\end{Prop}
\begin{proof}
For brevity, we will use the notations, definitions and basic results of
\cite[section 4]{erwanMaxEvenOvals}.
Suppose that there exists a curve $C$ contradicting the proposition. Without
loss of generality, we can assume that the
equation of $C$ is $y^2 + P(x)y + x=0$, where $P(x)$ is a real univariate
polynomial in 
$x$ of degree $2k-1$. The discriminant of $C$ seen as a polynomial in
$y$ is $R(x)=P^2(x) - 4x$, and it is clear that the topology of $C$ can be
recovered out of the root scheme realized by the polynomials $P^2(x)$,
$Q(x)=-4x$, 
and $R(x)$. One sees that $R(x)>0$ for $x\le0$, and since $C$ has an
$A_{2k-1}^\smbullet$ point, $R(x)$ must have a root of order
$2k$ close to which $R(x)$ is non positive. It follows that the polynomials
$P^2(x)$, $Q(x)$ and $R(x)$ realize the root scheme 
$$\begin{array}{ll}
\big((p,2b),(q,1),(r,a_1),(p,2b_1),(r,a_2),(p,2b_2),\ldots,(r,a_i),(p,2b_i),
(r,2k),\\
(p,2b_{i+1}),(r,a_{i+1}), (p,2b_k),(r,a_k)\big)
\end{array}$$
where $i$, $k$, $a_j$, $b$ and $b_j$ are some non negative integers, and $a_1>0$.
It is not
hard to see from the real rational graphs (or Dessins d'Enfants) that
this is equivalent to the existence of three real polynomials
$\widetilde P^2(x)$,
$\widetilde Q(x)$ and $\widetilde R(x)$ of
degree $4k-2$ and which realize the root scheme 
$$\begin{array}{ll}
\big((r,1),(p,2b),(r,a_1-1),(p,2b_1),(r,a_2),(p,2b_2),\ldots ,(r,a_i),(p,2b_i),
(r,2k), 
\\(p,2b_{i+1}),(r,a_{i+1}), (p,2b_k),(r,a_k)\big).
\end{array}$$ 
But then,
$\widetilde Q(x)=-\beta^2$ with $\beta$ a nonzero real number, and
$\widetilde R(x)=\widetilde P(x)^2 - \beta^2=(\widetilde
P(x)-\beta)(\widetilde P(x)+\beta)$. Now, the 
polynomials $\widetilde P(x)-\beta$ and $\widetilde P(x)+\beta$ are
relatively prime and 
of degree $2k-1$, so $\widetilde Q(x)$ cannot have a root of order $2k$.
\end{proof}

\end{Rem}

\bibliographystyle{amsalpha}
\bibliography{manyRealSings.bbl}

\end{document}